\documentclass[12pt]{amsart}
\usepackage{setspace}

\newtheorem{thm}{Theorem}[section]
\newtheorem{lem}[thm]{Lemma}
\newtheorem{prop}[thm]{Proposition}

\newtheorem{cor}[thm]{Corollary}
\newtheorem{rmk}[thm]{Remark}

\newenvironment{pf}{\begin{proof}}{\end{proof}}

\newcommand{\FZ}{\Bbb{Z}}  
\newcommand{\FQ}{\Bbb{Q}}  

\newcommand{\Q}{\Bbb{Q}}
\newcommand{\Z}{\Bbb{Z}}

\newcommand{\fQ}{\frak{Q}}

\begin{document}

\title[caliber number]{Caliber number of real quadratic fields}

\author{Byungheup Jun}
\author{Jungyun Lee }
\email{byungheup@gmail.com, lee9311@snu.ac.kr}
\address{Department of Mathematics. Korea Advanced Institute of Science and Technology}
\footnote{"This work was supported by the SRC Program of Korea Science and Engineering Foundation (KOSEF) grant funded by the Korea government(MEST) R11-2007-035-01001-0."}

\begin{abstract}
We obtain lower bound of caliber number of real quadratic field $K=\FQ(\sqrt{d})$ using splitting primes in $K$.  
We find all real quadratic fields of caliber number $1$ and find all real quadratic fields of caliber number $2$ if $d$ is not $5$ modulo $8$. In both cases,  we don't rely on 
the assumption on $\zeta_K(1/2)$. 
\end{abstract}
\maketitle
\tableofcontents
\section{Introduction} 


In \cite{Gauss}, Gauss had conjectured that there exist exactly nine  imaginary quadratic fields of class number $1$. 
Later, this was solved after diverse works of Stark, Heegner and Baker.

Further in this direction Goldfeld found an explicit lower bound of 
the class number of a given discriminant assuming existence of an elliptic curve on each imaginary quadratic field $K=\FQ(\sqrt{d})$ whose Hasse-Weil $L$-function have order of vanishing $3$ at $s=1$(cf. \cite{Goldfeld}).
Together with  Gross-Zagier's formula for $L'_K(E,1)$, Goldfeld's bound yields an explicit upper bound of a discriminant $|d|$ with $h(d)=h_0$, where $h(d)$ is a class number of $K=\FQ(\sqrt{d})$.
Finally, this gives an effective way of finding all imaginary quadratic fields with a given class number $h_0$.  

Contrary to imaginary quadratic case, in real quadratic field case, the same question remains still unanswered. It is believed that there are infinitely 
many real quadratic fields of class number $1$. 
As the first step has not been answered, at this moment, it does not make much sense to ask a similar generalization as above due to Stark, Heegner, Baker, Goldfeld, Gross and Zagier et al. 

If we replace  class number with caliber number, there is a room for a 
parallel generalization for real quadratic fields as in imaginary quadratic fields. 
Let $d$ be a positive square free integer and $D$ be a discriminant of 
the real quadratic field $K=\FQ(\sqrt{d})$.
Denote by $[A,B,C]$  a binary quadratic form $Q(X,Y)=AX^2+BXY+CY^2\in \FZ[X,Y]$. 
Then $GL_2(\FZ)$ acts on the set 
$\fQ(D)$ of primitive binary quadratic forms $[A,B,C]$ of discriminant $D=B^2-4AC$ by $S\circ Q(X,Y)=Q(aX+bY,cX+dY)$ for 
$S= \left(\begin{array}{cc}a & b \\c & d\end{array}\right)\in GL_2(\FZ)$ and $Q(X,Y)\in \fQ(D).$
Let $H(D)$ be the set of equivalent classes $\fQ(D)/GL_2(\FZ)$. 
The cardinality of $H(D)$ is the
class number $h(D)=h(d)$ of $\Q(\sqrt{d})$.

A quadratic form $[A,B,C]$ of discriminant $D$ is called {\it reduced} if the coefficients satisfy the following inequalities:
\begin{equation}
A>0,\,\,B<0,\,\,C<0,\,\,|B|<\sqrt{D},\,\,\sqrt{D}-|B|<2A<\sqrt{D}+|B|.
\end{equation}
Let $\fQ_{red}(D)$ be the set of reduced forms.  
The caliber number $\kappa(D)=\kappa(d)$ of $\Q(\sqrt{d})$ is the cardinality of  $\fQ_{red}(D)$.

For a real quadratic irrationality $w$ in $K$, its caliber $m(w)$ 
is simply the length of the periodic part in the continued fraction expansion.  
Let $\omega_Q$ be a root of $Q(X,1)$, then $m([Q])$ 
is actually a class invariant in such a way that $m([Q])=m(w_{Q_1})=m(w_{Q_2})$ if $Q_1$ and $Q_2$ are in the same class $[Q]\in H(D).$ 
And it is well-known that $m([Q])$ is the number of reduced forms in the class $[Q]$.
Thus the caliber number $\kappa(D)$ is  rewritten as follows:
$$
\kappa(D)= \sum_{[Q]\in H[D]} m([Q])
$$

In \cite{Gilles}, Lachaud obtained an effective lower bound of $\kappa(d)$ 
assuming $\zeta_K(\frac{1}{2}) \le 0$:
\begin{equation}
\kappa(d) > \frac{1}{8.46}\log(d-3).
\end{equation}
This is a real analogue of Goldfeld's work. Moreover, Lachaud determined all real quadratic fields with caliber number $1$ with assumption of $\zeta_K(\frac{1}{2})\le 0$.

We explain the content of this article.

In Section \ref{section:lowerbound}, 
we recall the definition of a set $S_D(A)$ and its cardinality $\rho_D(A)$ for 
a positive integer $A$. 
We write 
a lower bound and an upper bound of $\kappa(d)$ in terms $\rho_D(A)$.
Some further properties of $\rho_D(-)$ are studied  
to give a lower bound of the caliber number of a real quadratic field in terms of $D$ and  a rational prime that splits above. 
The other results of this paper relies on this estimate:\\

\noindent\textbf{Theorem \ref{lower}.}
\textit{Let $d$ be a positive square free integer and $K=\FQ(\sqrt{d})$ be a real quadratic field of discriminant $D$. 
Suppose  a rational prime $p$  splits in $K$.  Then 
$$\kappa(d)>2\Big{[}\frac{\log\frac{\sqrt{D}}{2}}{\log p}\Big{]}.$$}
  
In Section \ref{section:Determinantion}, we investigate 
the caliber number problem of real quadratic fields without the assumption on $\zeta_K(1/2)$.\\

\noindent\textbf{Theorem \ref{caliber1}.}
($\kappa(d)=1$) {\it
Suppose $d$ is a positive square free integer.
Then $\kappa(d)=1$ if and only if $d$ is one of the following: $2, 13, 29, 53, 173, 293$.}\\

Since we have related the lower bound of caliber number with a splitting prime and 
the values of $\rho_D(-)$, we further obtain an existence of splitting prime smaller
than $\sqrt{D}$ in case of $\kappa(d)\ne 1$.
We further study the $\kappa(d)=2$ problem for some cases.
We apply some results on class number problems of Richaud-Degert type 
due to Bir\'o , Byeon and the second named author (cf. \cite{Biro1}, \cite{Biro2}, \cite{Byeon-Kim}, \cite{Lee}), 
we list all real quadratic fields of caliber number one and all real quadratic fields $K=\FQ(\sqrt{d})$ 
when $d\not\equiv5$ modulo $8$ with caliber number two.\\

\noindent\textbf{Theorem \ref{caliber2}.}
($\kappa(d)=2$ with $d\not\equiv 5~(8)$){\it 
Real quadratic fields $K=\FQ(\sqrt{d})$ with $d\not\equiv5$ modulo $8$ with caliber number $2$ are the followings:
$$3, 6, 11, 38, 83, 227.$$}

\textit{Acknowledgment.} The authors thank to Sey Yoon Kim for introducing 
Gauss' work on classification of the reduced forms in real quadratic fields and for
useful discussions and ideas.


\section{Lower bound of caliber number}\label{section:lowerbound}

Throughout this article, $D>0$ denotes the discriminant of the real quadratic field $K=\Q(\sqrt{d})$.

For each positive integer $A$, we associate a set 
$$
S_D(A) := \{ [B] \in \Z/2A\Z | B^2 \equiv D \pmod{4A} \}
$$ 
and let $\rho_D(A)$ be the cardinality of 
$S_D(A)$.  

\begin{lem}\label{lemma}
Suppose $A<\frac{\sqrt{D}}2$. 
Then for a given  $[B]\in S_D(A)$, there exists 
a unique pair of integers $(B,C)$ such that $B\in [B]$ and $[A,B,C]$ is reduced. 
\end{lem}
\begin{pf}
Let $B_0$ be any integer representative of $[B]$. Then 
%
for any integer $k$ one can find uniquely an integer $C(k)$ satisfying
$$
D=(B_0+2Ak)^2-4AC(k).
$$

Moreover, we have for a unique integer $k_0$,
$$
-\sqrt{D}< B_0+2Ak_0 <2A-\sqrt{D}.
$$

If we set $B=B_0+2Ak_0$ and $C=C(k_0)$, from $A<\sqrt{D}/2$, one
can check easily $B<0$ and $2A < \sqrt{D}-B$.
\end{pf}

\begin{thm}\label{ineq}
The caliber number $\kappa(D)$ of $\FQ(\sqrt{D})$ satisfies
$$
\sum_{A<\frac{\sqrt{D}}{2}}\rho_D(A)\leq \kappa(D)\leq\sum_{A<\sqrt{D}}\rho_D(A).$$
\end{thm}

\begin{pf}
The lower bound is 
immediate 
from Lemma \ref{lemma}.
If a primitive quadratic form $[A,B,C]$ is reduced then for $B<\sqrt{D}$,
$$2A<B+\sqrt{D}.$$ 
Thus  we obtain $$A<\sqrt{D}.$$
This yields the upper bound of $\kappa(d)$.
\end{pf}


\begin{lem}\label{mollin}
 Let $d$ be a positive square free integer and $K=\FQ(\sqrt{d})$ with discriminant $D$. Set 
$$
\omega_D=\begin{cases}
     \frac{\sqrt{D}}{2} & D\equiv0 \pmod{4} \\
   \frac{1+\sqrt{D}}{2}   &D\equiv1 \pmod{4}.
\end{cases}
$$
Then an integral ideal  is of the form $[a,b+c\omega_D]$ for some positive integers $a,b,c$  such that $c |b, c|a$ and $ac | N(b+c\omega_D)$.
\end{lem}
\begin{pf}
See Thorem 1.2.1 and Definition 1.2.1 in \cite{Mollin}.
\end{pf} 

We say that an integral ideal $[A,\frac{B+\sqrt{D}}{2}]$ of $K$ is {\it primitive} 
if the integers $A, B$ satisfy 
$$B^2\equiv D \pmod{4A}.$$
(See Theorem 1.2.1 and Definition 1.2.1 in \cite{Mollin}.)

\begin{lem}\label{lee} 
Let $d$ be a positive square free integer and $K=\FQ(\sqrt{d})$ with discriminant $D$.\\
1) $\rho_D(A)$ equals the number of primitive ideals of $K$ with norm $A$.\\
2) Any integral ideal  can be written as $[fA,\frac{fB+f\sqrt{D}}{2}]$ for a positive integer $f$ and a primitive ideal $[A,\frac{B+\sqrt{D}}{2}].$
\end{lem}
\begin{pf}
1) Note that 
$$[A,\frac{B'+\sqrt{D}}{2}]=[A,\frac{B+\sqrt{D}}{2}]$$ 
if and only if  
$$B'\equiv B \pmod{2A}.$$ 
Thus for a primitive ideal $I$ of $K$, there exists exactly a unique pair of integers 
$(A,B)$ with $[B]\in S_D(A)$ 
and $I=[A,\frac{B+\sqrt{D}}{2}]$. 
If $I=[A,\frac{B+\sqrt{D}}{2}]$ is a primitive ideal,  
$$N(I)=A.$$
This completes the proof.\\
2) It is an immediate consequence of Lemma \ref{mollin}.  

\end{pf}

\begin{prop}\label{rho}
1) If $(n,m)=1$, $\rho_D(nm)=\rho_D(n)\rho_D(m)$. \\
2) 
For $p\not| D$,
$$\rho_D(p^{\alpha})=1+\chi_D(p),
$$
where $\chi_D$ be the Kronecker character(ie.  $\chi_D(\cdot)=(\frac{D}{\cdot})$).
For $p|D$,
$$\rho_D(p^{\alpha})=\begin{cases}
0, \quad\alpha>1\\
1, \quad\alpha=1 
\end{cases}$$
\end{prop}
\begin{pf}
1) It is clear from  Lemma 3.2 of pp. 48 in \cite{cox}.\\
2) From Lemma \ref{lee}, we find that 
$$\sum_{n}\rho_D(n)n^{-s}=\zeta(2s)^{-1}\zeta_K(s).$$ 
The Euler factor at $p$  of $\sum_{n}\rho_D(n)n^{-s}$ is 
$$
1+\sum_{n=1}^{\infty}\frac{\rho_D(p^n)}{p^{ns}}.
$$
If $\chi_D(p)=1$ (resp. $\chi_D(p)=-1$, $\chi_D(p)=0$)  then the Euler factor at $p$ 
of $\zeta(2s)^{-1}\zeta_K(s)$ is
$$1+\sum_{n=1}^{\infty}\frac{2}{p^{ns}}\quad (\text{resp.} ~~1,  1+p^{-s}).$$
By comparing  the Euler factors of $\sum_{A}\rho_D(A)A^{-s}$ and $\zeta(2s)^{-1}\zeta_K(s)$, we can prove Proposition.
\end{pf}

\begin{thm}\label{lower}
Let $d$ be a positive square free integer and $K=\FQ(\sqrt{d})$ be a real quadratic field of discriminant $D$.
Suppose a rational prime $p$ splits in $K$. 
Then
$$\kappa(d)>2\Big{[}\frac{\log\frac{\sqrt{D}}{2}}{\log p}\Big{]}.$$
\end{thm}
\begin{pf}
It suffices to show the theorem for the smallest splitting prime.
Let $p_1$ be the smallest prime that splits in $K$. 

From Theorem \ref{ineq}, we have:
$$
\sum_{p_1^{\alpha}<\frac{\sqrt{D}}{2}}\rho_D(p_1^{\alpha})\leq\sum_{A<\frac{\sqrt{D}}2}\rho_D(A)<\kappa(d).
$$
Lemma \ref{rho} implies that $\rho_D(p_1^\alpha)=0$ for any $\alpha$. 

Therefore,
\begin{equation*}
\begin{split}
\sum_{p_1^{\alpha}<\frac{\sqrt{D}}2} \rho_D(p_1^\alpha) &= 2\cdot 
\text{(the number of $\alpha$'s:$p_1^\alpha < \sqrt{D}/2$ })\\
&= 2 \Big{[}\frac{\log{\frac{\sqrt{D}}{2}}}{\log p_1}\Big{]}.
\end{split}
\end{equation*}
This completes the proof.
\end{pf}

\begin{cor}
Suppose $d\equiv 1 \pmod{8}$ be a positive square free integer. 
Then 
$$2^{\kappa(d)+4}>d.$$
\end{cor}

\begin{rmk}
The result of this section is comparable to 
Section 22.5 of \cite{Iwaniec}. 

For an imaginary quadratic field of discriminant $D<0$,  one has
\begin{equation}\label{classnumber}
\sum_{A\leq\sqrt{\frac{|D|}{4}}}\rho_D(A)\leq h(D)\leq\sum_{A\leq\sqrt{\frac{|D|}{3}}}\rho_D(A).
\end{equation}
As $h(D)=\kappa(D)$ for $D<0$, the above extends Proposition \ref{rho} to 
negative discriminant case. 
The inequality $(\ref{classnumber})$ turns out to a lower bound of $h(D)$ in terms of $\log p_r$
$$ \log p_r \geq \frac{\log \frac{|D|}{4}}{\sqrt[r]{2 h(D) r!}}$$
where $p_1<p_2<\cdots <p_r$ are the first $r$ splitting primes.   

\end{rmk}

\section{Determination of real quadratic fields with small caliber numbers}\label{section:Determinantion}

In this section, we determine all the real quadratic fields $K=\FQ(\sqrt{d})$ with caliber  $1$ and  $\FQ(\sqrt{d})$ with caliber number $2$ when $d\not\equiv 5$ modulo $8$. For the determination, we don't assume $\zeta_K(1/2)\le 0$.

For both cases, we need some ideas on continued fractions of quadratic irrationalities.
For general and precise idea on continued fractions we refer the readers to
\cite{Mollin}, \cite{Silverman}, \cite{vdGeer}, etc.

Consider a real quadratic irrationality $x$.
The caliber $m(x)$ is simply the length of the periodic part 
in the continued fraction expansion.  $x$ is said to be reduced 
if $x>1$ and $-1<x'<0$. It is well known that the reduced elements 
$x$ has purely periodic continued fraction expansion.

\begin{prop}\label{prop1}
Let $d$ be a positive square free integer. 
Then $K=\FQ(\sqrt{d})$ is of caliber one only if it has class number one and
$$d=n^2+4\,\,\text{or}\,\, n^2+1.$$
\end{prop}
\begin{pf}
Suppose $K=\Q(\sqrt{d})$ is of caliber $1$.
Since $h(d)\le \kappa(d)$, clearly $h(d)=1$.

Suppose now $x$ is reduced with period $1$ then for a positive integer $r$,
$x$ satisfies
$$x=r+\frac{1}{x}.$$
Solving the above equality, we get
$$x=\frac{r+\sqrt{r^2+4}}{2}.$$

Suppose $d\equiv 1 \pmod{4}$. 
Since $\kappa(d)=h(d)=1$, $[1,\frac{1+\sqrt{d}}{2}]$ is  principal and $m(\frac{1+\sqrt{d}}{2})=1$. 
As $(\frac{1+\sqrt{d}}{2}-[\frac{1+\sqrt{d}}{2}])^{-1}$ is reduced, 
we have for a positive integer $r$
$$
(\frac{1+\sqrt{d}}2 - \Big[\frac{1+\sqrt{d}}{2}\Big])^{-1} = \frac{r+\sqrt{r^2+4}}2.
$$
Thus
$$\frac{1+\sqrt{d}}{2}-\Big[\frac{1+\sqrt{d}}{2}\Big]=\frac{-r+\sqrt{r^2+4}}{2}.$$ 
From above equation, we have 
$$d=\Big{(}2\Big{[}\frac{1+\sqrt{d}}{2}\Big{]}-1\Big{)}^2+4.$$

If $K=\FQ(\sqrt{d})$ with $d\equiv 2,3 $ modulo $4$ of caliber number $1$, 
similarly we can conclude that $d$ is of the form $n^2+1$ for an integer $n$.
\end{pf}

In the above case of $d$ (i.e.. $d=n^2+4$ or $n^2+1$), Bir\'o found the full list of $d$ with $h(d)=1$. 

\begin{prop}[Bir\'o \cite{Biro1}, \cite{Biro2}]
\label{prop2}\item{I.} Let $d=n^2+4$ be a square free integer. $h(d)=1$  if and only if 
\begin{equation}\label{list1}
d=13,19,53,173,293.
\end{equation} 
\item{II.} Let $d=n^2+1$ be a square free integer. $h(d)=1$ if and only if 
\begin{equation}\label{list2}
d=2,17,37,101,197,677.
\end{equation}
\end{prop}

Combining Proposition \ref{prop1} and Proposition \ref{prop2}, we can 
list all $d$ with $\kappa(d)=1$:

\begin{thm}[$\kappa(d)=1$] \label{caliber1}
Suppose $d$ is a positive square free integer.
Then $\kappa(d)=1$ if and only if $d$ is one of the following: $2, 13, 29, 53, 173, 293$.
\end{thm}

\begin{cor} 
Suppose that $d$ is a rational prime that is not in
$$S:=\{2,13,29, 53, 173,293\}$$.
Let $D$ be the discriminant of $\FQ(\sqrt{d})$. 
Then there exists 
a rational prime that splits in $\Q(\sqrt{d})$ and smaller than or equal to $\sqrt{D}.$   
\end{cor}
\begin{pf}
After Theorem \ref{caliber1}, we know that $\kappa(d)\geq 2$ for $d\not\in S$.  

Let $A\ne 1$ be a positive integer smaller than $\sqrt{D}$.
If we suppose conversely that there no prime $\leq \sqrt{D}$ splits above,
then  $\rho_D(A)=0$ or $1$. 
If the multiplicity of a prime factor $p$ of $A$ is greater than $2$, $\rho_D(A)=0$.
And if a prime factor $p$ of $A$ inerts in $\Q(\sqrt{d})$, $\rho_D(A)=0$ from 
the multiplicative property of $\rho_D$ as in Proposition \ref{rho}.

Since $S_D(1)\subseteq \Z/2\Z$ cannot contain $0$, $\rho_D(1) <2$. 
As $D$ is either $d$ or $4d$, $d > \sqrt{D}$.

Therefore, 
we have   
$$\kappa(d) \le \sum_{A\leq\sqrt{D}}\rho_D(A) \le \rho_D(1) \leq1.$$
This contradicts to our assumption. 
\end{pf}


Now we move to the case of $\kappa(d)=2$.

A positive square free integer $d=n^2+r$ with $r|4n$ is said to be of Richaud-Degert type.
Lemma \ref{lem44} and Proposition \ref{prop45}
imply that if 
$\kappa(d)=2$ then $d$ is necessarily of Richaud-Degert type:

\begin{lem}\label{lem44}
Let $d=n^2+1$ be a square free integer. If $d\equiv 2 \pmod{4}$ then the ideal 
$[2,\sqrt{d}]$ is not principal ideal except $d=2$ and if $d\equiv 3$ modulo $4$ then the ideal 
$[2,1+\sqrt{d}]$ is not principal ideal
if  $d\equiv 1$ modulo $8$ then the ideal 
$[2,\frac{1+\sqrt{d}}{2}]$ is not principal ideal except $d=17$.
\end{lem}
\begin{pf}
See the proof of Theorem 2.6 in \cite{Byeon-Kim}.
\end{pf}

\begin{prop}\label{prop45}
Let $d\not\equiv 5$ modulo $8$ be a positive square free integer. The field $\FQ(\sqrt{d})$ is of caliber number 2 then  $d$ is a Richaud-Degert type with class number one.
\end{prop}

\begin{pf}
Suppose $K=\FQ(\sqrt{d})$ has caliber number $2$ and class number $2$. 
Then the principal ideal has caliber $1$. 
Thus from the proof of Proposition \ref{prop1}, 
we know that $d$ is either $n^2+1$ or $n^2+4$.
As we assumed that 
$d\not\equiv 5 \pmod{8}$, we can exclude  the case $d=n^2+4$.
From Lemma \ref{lem44}, if $n^2+1\equiv 2 \pmod{4}$,
$[1,\sqrt{n^2+1}]$ and $[2,\sqrt{n^2+1}]$ represent two distinct ideal classes of caliber $1$. 
Thus $[2,\sqrt{n^2+1}]\sim [1,\sqrt{n^2+1}/2] \sim [1,x]$, where 
$$x^{-1}=\frac{\sqrt{n^2+1}}{2}-\Big[\frac{\sqrt{n^2+1}}{2}\Big]=\frac{-r+\sqrt{r^2+4}}{2}$$
for a positive integer $r$.
In the above, comparing the rational parts, one can see that 
\begin{eqnarray*}
&r=2\Big{[}\frac{\sqrt{n^2+1}}{2}\Big{]},\\
&n^2+1=r^2+4.
\end{eqnarray*}
This contradict to the assumption that $d=n^2+1$ is square free. 
Similarly, for the rest cases $d\equiv 3 \pmod{4}$ or $d\equiv 1\pmod{8}$, we obtain contradiction.

Therefore, if $d\not\equiv 5\pmod{8}$ and $\kappa(d)=2$, then $h(d)=1$. 

Suppose now $x$ is a reduced quadratic irrationality of caliber $2$. Then 
\begin{equation}
x = a+\frac{1}{b+\frac{1}{x}}.
\end{equation}
for two distinct positive integers $a$ and $b$.
Solving the above equation, we 
obtain 
\begin{equation}\label{cal2}
x=\frac{ab+\sqrt{a^2b^2+4ab}}{2b}.
\end{equation}

Consider the case $d\equiv 1\pmod{8}$. 
Since
$\kappa(d)=2$ implies $h(d)=1$, 
$[1,\frac{1+\sqrt{d}}{2}]$ is principal and  $m(\frac{1+\sqrt{d}}{2})=2.$
Thus from the equation (\ref{cal2}), we have
\begin{equation}\label{caliber2}
x^{-1}= \frac{1+\sqrt{d}}{2}-\Big[\frac{1+\sqrt{d}}{2}\Big]=\frac{-ab+\sqrt{a^2b^2+4ab}}{2a}.
\end{equation} 
And (\ref{caliber2}) implies that 
\begin{eqnarray*}
&b=2\Big{[}\frac{1+\sqrt{d}}{2}\Big{]}-1,\\
&d=b^2+4\frac{b}{a}.
\end{eqnarray*}
Thus we find that $d$ is of the form $n^2+r$ with $r|4n$.

Similarly, for a square free integer $d\equiv 2,3$ modulo $4$, 
we conclude that if $\kappa(d)=2$, then 
$d$ is of the form $n^2+r$ with $r|2n$.
\end{pf}

For Richaud-Degert types of $d\not\equiv 5 \pmod{8}$, 
we recall a class number $1$ criterion by Byeon and Kim (cf. \cite{Byeon-Kim}):
\begin{prop}[Byeon-Kim]\label{B-K}
Let $K=\FQ(\sqrt{d})$ be a real quadratic field of R-D type and $h(d)$ be the class number of $K$. Then
\begin{itemize}
\item[I.]
$d=n^2+r\equiv 2, 3$ (mod 4)
\begin{itemize}
\item[(i)]
$|r|\neq 1, 4,$ $h(d)>1$  except $r=\pm 2$
\item[(ii)]$|r|=1,$ $h(d)>1$ except $d=2,3$
\end{itemize}
\item[II.]
$d=n^2+r\equiv 1$ (mod 8)
\begin{itemize}
\item[(i)] $|r|\neq 1, 4$ $h(d)>1$  except $d=33$
\item[(ii)]
$|r|=1$ (hence $r=1$ and $n$ even)
$h(d)>1$ except $d=17$.
\end{itemize}
\end{itemize}
\end{prop}

After the Proposition \ref{B-K}, if $d$ is a Richaud-Degert type of $h(d)=1$, then 
$d= n^2\pm 2, 2,3,33$ or $17$.

For Richaud-Degert type of $n^2\pm2$,
the second named author found the whole list of $d$ with $h(d)=1$ in \cite{Lee}: 
\begin{prop}[Lee]\label{Lee}
Let $d=n^2\pm2$ be a square free integer. Then $h(d)=1$ if and only if 
$$d=3, 6, 7, 11, 14, 23, 38, 47, 62, 83, 167, 227, 398. $$ 
\end{prop}


Finally,  we obtain the following list of $d(\not\equiv5\pmod{8})$ with $\kappa(d)=2$:
\begin{thm}[$\kappa(d)=2$ with $d\not\equiv 5~(8)$]\label{caliber2}
Real quadratic fields $K=\FQ(\sqrt{d})$ with $d\not\equiv5$ modulo $8$ with caliber number $2$ are the followings:
$$3, 6, 11, 38, 83, 227.$$
\end{thm}



\begin{thebibliography}{99}


\bibitem{Biro1} A. Bir\'o,
{\it Yokoi's conjecture}, Acta Arith. {\bf106} (2003), 85-104.


\bibitem{Biro2} A. Bir\'o,
{\it Chowla's conjecture}, Acta Arith. {\bf107} (2003), 179-194.


\bibitem{Byeon-Kim} D. Byeon and H. Kim,
{\it Class number 2 Criteria for real quadratic fields of Richaud-Degert type}, 
Journal of Number Theory {\bf 62} No {\bf 2} (1997) 257-272.

\bibitem{Goldfeld} D. Goldfeld,
{\it The class number of quadratic fields and the conjectures of Birch and Swinnerton-Dyer}, Ann. Scuola Norm. Sup. Pisa Cl. Sci. (4) {\bf3} (1976), 624-663.

\bibitem{Gauss} C. F. Gauss, 
{\it Disquisitiones arithmeticae},
 Translated by Arthur A. Clarke, Springer-Verlag, New York, 1986.

\bibitem{cox} David A. Cox,
{\it Primes of the Form $x^2+ny^2$}, Pure and Applied Mathematics.

\bibitem{Gilles} G. Lachaud,
{\it On real quadratic fields,} Bulletin of the A.M.S. Volume {\bf17}, Number {\bf2}, (1987), 307-311.



\bibitem{Lee} J. Lee,
{\it The complete determination of wide Richaud-Degert type which is not 5 modulo 8 with class number one}, to appear in Acta Arith.
\bibitem{Mollin} Mollin,
{\it Quadratics}, CRC Press Series on Discrete Mathematics and its Applications. CRC Press, Boca Raton, FL, 1996. xx+387 pp. 

\bibitem{Iwaniec} Henryk Iwaniec,
{\it Analytic Number Theory}, American Mathematical Society.

\bibitem{Siegel} C. L. Siegel,
{\it \"Uber die Classenzahl quadratischer Zahlk\"orper,} Acta Arith. {\bf 1.} (1934), 83-86.

\bibitem{Silverman} J. Silverman,
{\it A friendly introduction to number theory}, 3rd ed. Pearson Prentice Hall (2006), 439 pages. 

\bibitem{vdGeer} G. van der Geer, {\it Hilbert modular surfaces}....
 Ergebnisse der Mathematik und ihrer Grenzgebiete (3), {\bf 16}. Springer-Verlag, Berlin, (1988). 
\bibitem{Vijaya} T. Vijayayaghavan,
{\it Periodic simple continued fractions,} Proc London Math. Soc. {\bf 26} (1927), 403-414.

\end{thebibliography}
\end{document}